\documentclass[12pt,a4paper]{article}
\usepackage[latin1]{inputenc}
\usepackage{amsmath}
\usepackage{amsfonts}
\usepackage{amssymb}
\usepackage{array}

\usepackage[notcite,notref]{showkeys}

\newcommand{\mb}{ \bar{m} }

\newcommand{\R}{  \mathbb{R} }

\newcommand{\ZZ}{\mathbb{Z} }
\newcommand{\Z}{\mathbb{Z} }
\newcommand{\h}{  {h}  }
\newcommand{\vp}{ \varphi  }
\newcommand{\sG}{ {\mathcal G} }
\newcommand{\sX}{ {\mathcal X} }
\newcommand{\sU}{ {\mathcal U} }
\newcommand{\sF}{ {\mathcal F} }

\newcommand{\half}{ \frac 1 2  }

\newtheorem{proposition}{Proposition}
\newtheorem{theorem}{Theorem}[section]

\newtheorem{definition}[theorem]{Definition}

\numberwithin{equation}{section}

\author{Moritz Duembgen \thanks{University of Cambridge, 
Statistical Laboratory,
Wilberforce Road, Cambridge, CB3 0WB, UK} \\ University of Cambridge
\and L.~C.~G. Rogers \thanks{
Corresponding author.
University of Cambridge, 
Statistical Laboratory,
Wilberforce Road, Cambridge, CB3 0WB, UK} \\ University of Cambridge}
\title{{\bf The joint law of the extrema, final
value and signature of a stopped random walk.}  }
\begin{document}
\maketitle

\begin{abstract}
A complete characterization of the possible joint distributions
of the maximum and terminal value of  uniformly integrable
martingale has been known for some time, and the aim of this 
paper is to establish a similar characterization for continuous
martingales of the joint law of the minimum, final value, and 
maximum, along with the direction of the final excursion. We 
solve this problem completely for the discrete analogue, that of 
a simple symmetric random walk stopped at some almost-surely
finite stopping time. This characterization leads to robust hedging
strategies for derivatives whose value depends on the maximum,
minimum and final values of the underlying asset.
\end{abstract}

\section{Introduction.}\label{intro}
Suppose given $\h >0$, and 
suppose that $(\xi_t,\sF_t)_{t\in \h\Z^+} $ is a symmetric
simple random walk on the grid $\h\Z$, started at zero.
Define $S_t \equiv \sup_{s\leq t} \xi_s$,  
$I_t \equiv \inf_{s\leq t} \xi_s$,
 $g^+_t \equiv \inf\{ u \leq t: \xi_u = S_u\}$, 
$g^-_t \equiv \inf\{ u \leq t: \xi_u = I_u\}$, and let
\begin{eqnarray}
\sigma_t &=& +1 \qquad \hbox{\rm if $g^+_t > g^-_t$}
\nonumber
\\
&=& -1 \qquad \hbox{\rm else}.
\label{sigdef}
\end{eqnarray}
The process $S$ records the running maximum of the martingale, and 
the process $\sigma$ records whether the martingale is currently on
an excursion down from its running maximum ($\sigma = +1$) or on an
excursion up from its running minimum ($\sigma = -1$).  We refer
to the process $\sigma$ as the {\em signature} of the random walk.

Suppose that $T$ is an almost-surely finite $(\sF_t)$-stopping time,
 and write
 \[
        X_t \equiv \xi_{t\wedge T}
 \]
 for the stopped process. The paper
is concerned with the possible joint laws $m$ of the quadruple
$(I_T, X_T, S_T, \sigma_T)$, which we will abbreviate to 
$(I,X,S,\sigma)$ where no confusion may arise.

 Clearly the
law $m$ must be defined on the set $
\sX \equiv (-h\Z^+) \times \h\Z \times
\h\Z^+ \times\{ -1,+1\}$, 
and evidently  we must have
$m( I \leq X \leq S) = 1$;
but beyond this, 
is it possible to state a set of {\em necessary and sufficient
conditions} for a probability $m$ on $\sX$ to be the joint distribution
of $(I_t,X_T,S_T,\sigma_T)$? The motivation
for this attempt is twofold. Firstly, the joint law of $(X,S)$ has been
characterized completely (for general local martingales, 
not assumed to be continuous or uniformly integrable)
 in \cite{R1}; can the methods of that
paper be extended to deal with the running minimum also? The second
reason to look at this problem is the interesting recent work of 
Cox \& Obloj \cite{CO} which finds extremal martingales for various
derivatives whose payoffs depend on the maximum, minimum and terminal
value of the underlying asset. This builds to some extent on the
earlier work of Hobson and others 
(\cite{Hobson1}, \cite{BHR1}, \cite{BHR2}), 
which addresses similar questions
for derivatives whose payoffs depend only on the maximum and terminal 
value of the underlying asset. Many of the results of this literature
can be derived alternatively using the results of \cite{R1}, 
by converting the problem into a linear program. This approach is 
more general, but leads to less explicit answers in the specific
instances analyzed to date.

What we shall find here is that it is possible to generalize the 
results of \cite{R1} to cover the joint law of $(I,X,S,\sigma)$, 
but that the statements are more involved. For this reason, 
we shall restrict our analysis to  a  symmetric simple random walk 
taking values in a grid $\h\ZZ$ for some $\h>0$, stopped at an
almost-surely finite stopping time.  The main result is presented
in Section \ref{S1}. The proof of necessity is in Section \ref{nec},
and requires only the judicious use of the Optional Sampling
Theorem. The proof of sufficiency, in Section \ref{suff}, is
constructive, and requires suitable modification of some of the 
techniques of \cite{R1}. We then show in Section \ref{hedge} how
this characterization can lead to robust hedging schemes and
extremal prices for derivatives whose payoff depends on the
maximum, minimum, terminal value and signature.

\section{The main result.}\label{S1}
We take a symmetric simple random walk $(\xi_t,\sF_t)_{t\in \h\Z^+} $ 
on $\h\Z$ for some fixed $h>0$; in general, the filtration
$(\sF_t)$ is larger than the filtration of the random walk, to 
allow for additional randomization. Stopping $\xi$ at the 
almost-surely finite stopping time $T$ creates the martingale
$X_t = \xi_{t \wedge T}$.
 We use the notation of the Introduction,
and notice that 
\begin{equation}
g^+_t \equiv \sup\{ u \leq t: S_u > S_{u-\h}\},
\qquad g^-_t \equiv \sup\{ u \leq t: I_u < I_{u-\h}\},
\label{eq2}
\end{equation}
emphasizing the fact that we are dealing with {\em strict}
ascending/descending ladder epochs, to use the language of 
Feller \cite{Feller}.  The process $\sigma$ is defined as
before at \eqref{sigdef}.

\begin{definition}
We say that the probability measure $m$ on 
$\sX \equiv -\h\ZZ^+ \times \h\ZZ \times \h\ZZ^+ \times \{-1,+1\}$ is
{\em consistent} if there is some 
almost-surely finite  $(\sF_t)$-stopping time $T$ 
such that $m$ is the 
law of $(I_T,X_T,S_T,\sigma_T)$.
\end{definition}

    \subsection{Necessity.}\label{nec}
For $x \in \h\ZZ$ we define the hitting time
\begin{equation}
   H_x = \inf\{ u: \xi_u = x\},
   \label{Hdef}
\end{equation}
with the usual convention that the infimum of the empty set is 
$+\infty$. In what follows, we will let $a$, $b$ stand for 
two generic members of $\h \Z^+$, and will be studying the exit
time $H_b \wedge H_{-a} \equiv \inf\{u: \xi_u \notin (-a,b)\}$
and related stopping times.  The measure $m$ says nothing
directly about these stopping times, but by way of the 
Optional Sampling Theorem we are able to deduce quite a lot
of information about them if the law $m$ is consistent.
Indeed, assuming that $m$ is consistent, we are able to 
find the probability that $H_{-a}<H_b$ (for example) in terms of 
$m$-expectations of functions defined on $\sX$.  The expressions
derived make perfectly good sense even if $m$ is not consistent, 
but it may be that the expressions do not in general satisfy 
positivity or other properties which would hold if $m$ were
consistent.  For this reason, we will denote by $\mb(Y)$
the expression for 
the $m$-expectation of a random variable $Y$ which would be 
correct if $m$ were consistent; if $m$ is not consistent, 
all we have is an algebraic expression without the desired
probabilistic meaning, and the use of the symbol $\mb$
warns us not to assume properties which need not hold.

The first result we need is the following, which illustrates
the use of this notational convention.

\begin{proposition}\label{prop1}
For any $a, b \in \h \Z^+$ we have
\begin{eqnarray}
 \mb(H_b < H_{-a})  &=&
\frac{ a - m(a+X; S<b, I>-a) }{a+b}\;
\equiv \vp(b,-a) ,
\label{m+}
\\
\mb(H_{-a} < H_b)&=&\frac{b - m(b-X; S<b, I>-a) }{a+b}\;
\equiv \vp(-a,b)  .
\label{m-}
\end{eqnarray}
\end{proposition}

\medskip\noindent
{\sc Proof.} We use the Optional Sampling
Theorem at the time $H_b \wedge H_{-a}$ to derive the two
equations
\begin{eqnarray}
1 &=& \mb(H_{-a} < H_b) + \mb(H_b < H_{-a}) + m(S<b, I>-a)
\label{p1}
\\
0 &=& -a\, \mb(H_{-a} < H_b) +b\, \mb(H_b < H_{-a}) + m(X;S<b, I>-a). 
\label{OST1}
\end{eqnarray}
Solving this pair of linear equations leads to the conclusion that
\begin{eqnarray*}
\mb(H_b < H_{-a}) &=& \bigl\lbrace\, a - m(a+X; S<b, I>-a) \,\bigr
\rbrace/(a+b) \;  ,
\\
\mb(H_{-a} < H_b) &=& \bigl\lbrace\, b - m(b-X; S<b, I>-a) \, \bigr
\rbrace/(a+b) \;   ,
\end{eqnarray*}
as claimed.

\hfill$\square$

\medskip\noindent
If $m$ is consistent, then we would have for any $a, b \in \h\Z^+$
not both zero that
\begin{eqnarray*}
\mb(H_{-a} < H_b < H_{-a-\h}) &=& \mb(H_{-a} \leq H_b < H_{-a-\h})
\\
&=& \mb(H_b < H_{-a-\h}) - \mb(H_b < H_{-a})
\\
&=& \mb(H_b < \infty, I(H_b) = -a).
\end{eqnarray*}
This is because on the event $\{ H_{-a} < H_b < H_{-a-\h} \} $
the hitting time $H_b$ is finite, and so cannot be equal to 
$H_{-a}$; the second equality follows from the inclusion
$\{ H_b < H_{-a} \} \subseteq \{H_b < H_{-a-\h} \}$. We will 
therefore introduce the notation
\begin{eqnarray}
\psi_+(-a,b) &=& \vp(b,-a-\h) - \vp(b,-a),
\label{psip}
\\
\psi_-(-a,b) &=& \vp(-a,b+\h) - \vp(-a,b).
\label{psim}
\end{eqnarray}
Notice that $\psi_+(-a,b)$ is {\em defined} 
as an algebraic expression in terms of $m$
via \eqref{psip} and \eqref{m+}; if $m$ is {\em consistent}, then 
$\psi_+(-a,b) $ is equal to  $ \mb(H_b<\infty, I(H_b) = -a)$,
but no such interpretation holds in general.

\vskip 0.2 in

The  necessary condition we derive comes from considering what
may happen if the event $B_+ = \{ H_b < \infty, I(H_b) = -a\}$
occurs. When this event occurs, the martingale $X$ does reach 
$b$ before being stopped, and at that time $H_b$ the minimum
value is $-a$. Thereafter, one of three things will happen:
\begin{itemize}
\item[(i)] $X$ reaches $b+\h$ before reaching $-a-\h$ and 
before $T$;
\item[(ii)] $T$ happens before $X$ reaches either $-a-\h$
or $b+\h$;
\item[(iii)] $X$ reaches $-a-\h$ before reaching $b+\h$ and
before $T$.
\end{itemize}
The next result derives a necessary condition from the Optional
Sampling Theorem applied  at $H_{-a-\h} \wedge H_{b+\h} \wedge T$.

\begin{proposition}\label{prop2}
Define the events
\begin{equation}
B_+ = \{ H_b < \infty, I(H_b) = -a)\}, \qquad
B_- = \{ H_{-a}<\infty, S(H_{-a}) = b  \},
\label{Bdef}
\end{equation}
set $p_\pm = \mb(B_\pm) = \psi_\pm(-a,b)$, 
and set
\begin{equation}
p_{+0} = m(S = b, I = -a, \sigma  =+1),
\qquad
p_{-0} = m(S = b, I = -a, \sigma  =-1).
\end{equation}
If we denote
\begin{equation}
v_{\pm} \equiv  \frac{m(X; S=b,I=-a,\sigma = \pm 1)}{p_{\pm 0}  }
\equiv m(X \, \vert \, S=b,I=-a,\sigma = \pm 1),
\label{vdef}
\end{equation}
then the conditions\footnote{If either of $p_\pm$ is zero, then 
the inequalities \eqref{nc5}, \eqref{nc6} have to be understood in 
cross-multiplied form, when they state vacuously that $0 \leq 0$.}
\begin{eqnarray}
  \frac{p_{+0}}{p_+} & \leq &\frac{\h}{b+\h-v_+}
   \label{nc5}
   \\
   \frac{p_{-0}}{p_-}& \leq &\frac{\h}{a+\h+v_-}
   \label{nc6}
\end{eqnarray}
are necessary for $m$ to be consistent.
\end{proposition}

\medskip\noindent
{\sc Proof.}
We introduce the notation 
\begin{gather}
p_{++} = \mb( H_{-a}< H_b<H_{b+\h} < H_{-a-\h} ),
\quad
p_{+-} = \mb( H_{-a}< H_b<H_{-a-\h} < H_{b+\h} ),
\nonumber
\\
p_{--} = \mb( H_{b}< H_{-a}<H_{-a-\h} < H_{b+\h} ),
\quad
p_{-+} = \mb( H_{b}< H_{-a}<H_{b+\h} < H_{-a-\h} ).
\nonumber
\end{gather}
Using the Optional Sampling Theorem,
we have similarly to \eqref{p1}, \eqref{OST1} the equations
\begin{eqnarray}
p_+ &=& p_{++} + p_{+0} + p_{+-}
\label{p2}
\\
bp_+ &=& (b+\h) p_{++}
-(a+\h) p_{+-}
+m(X; S=b,I=-a,\sigma = +1).
\label{OST2}
\end{eqnarray}
If we write $\tilde p_{xy} = p_{xy}/p_x$ for $x\in 
\{  -, + \}$,  $y \in \{ -, 0 ,  +\}$
the equations \eqref{p2}, \eqref{OST2} are expressed more simply
in conditional form:
\begin{eqnarray}
1 &=& \tilde p_{++} + \tilde p_{+-} + \tilde p_{+0}
\label{p3}
\\
b &=& (b+\h)\tilde  p_{++}
-(a+\h) \tilde p_{+-}
+\tilde p_{+0} v_+.
\label{OST3}
\end{eqnarray}
The value of $p_{+0}$ is known from $m$, as is the value
of $v_+$, and since we assume that $m$ is consistent the
values of $p_\pm = \psi_\pm(-a,b)$ are also known from $m$.
  Therefore we can solve the linear system
\eqref{p3}, \eqref{OST3} to discover 
\begin{eqnarray}
\tilde p_{++} &=& \frac{b+a+\h - (a+\h+v_+)\, \tilde p_{+0}}{b+a+2\h}
\label{p++}
\\
\tilde p_{+-} &=& \frac{h - (b+\h-v_+)\,  \tilde p_{+0}}{b+a+2\h}.
\label{p+-}
\end{eqnarray}
In order that $\tilde p_{+-}$ as given by \eqref{p+-} should be 
non-negative, we require that 
\begin{equation}
   \tilde p_{+0} \equiv \frac{m(S = b, I = -a, \sigma  =+1)}
   {p_+ } \leq \frac{\h}{b+\h-v_+},
   \label{p+0}
\end{equation}
which is condition \eqref{nc5}. Necessity of \eqref{nc6} is 
derived similarly.

\hfill$\square$

\medskip\noindent
{\sc Remarks.} (i) The necessary conditions \eqref{nc5}, \eqref{nc6}
come from the requirement that $\tilde p_{+-}$ and $\tilde p_{-+}$
should be non-negative. Do we know for sure that $\tilde p_{++}$
and $\tilde p_{--}$ are non-negative?  The definition \eqref{vdef}
of $v_\pm$ guarantees that $-a \leq v_\pm \leq b$, so if \eqref{p+0}
holds then we know that $\tilde p_{+0} \leq 1$.  From \eqref{p++}
we see then that $\tilde p_{++} \geq 0$. Since all the summands on the
right-hand side of \eqref{p3} are non-negative, we learn that
they are probabilities summing to 1.

\medskip\noindent
(ii) Notice that we have two expressions for
 $\mb(H_{b+\h}<\infty,
I(H_{b+\h}) = -a)$, either as $p_{++} + p_{-+}$, or as 
$\psi_+(-a,b+\h)$.
Confirming that these are the same is an important step in the proof 
of sufficiency.

\bigbreak
    \subsection{Sufficiency.}\label{suff}
We have now identified necessary conditions 
\eqref{nc5} and \eqref{nc6} for $m$ to be consistent.
The main result of this paper is that these conditions are also
sufficient.
\begin{theorem}\label{thm1}
The probability measure $m$ on $\sX \equiv
-\h\ZZ^+ \times \h\ZZ \times \h\ZZ^+ \times \{-1,+1\}$ is
consistent
if and only if  $ m(I \leq X\leq S) = 1$ and necessary conditions 
\eqref{nc5} and \eqref{nc6} hold.
\end{theorem}

\bigbreak\noindent
{\sc Proof.}  Necessity has been proved, so 
what remains is to show that conditions 
\eqref{nc5} and \eqref{nc6} are 
sufficient. Not surprisingly, the proof of this is constructive.

We require a probability space $(\Omega, \sF, P)$ rich enough to 
carry an IID sequence $U_0, U_1, \ldots$ of $U[0,1]$ random 
variables, and an independent standard Brownian motion $(B_t)$.
Let $\sU = \sigma(U_0, U_1, \ldots)$, and let $(\sG_t)$ be
 the usual augmentation of the  filtration
$(\sU \vee \sigma(B_s: s\leq t))$.  Define $(\sG_t)$-stopping times
\begin{equation*}
\alpha_0 \equiv 0, \qquad
\alpha_{n+1} \equiv \inf\{ t> \alpha_n : |B_t-B_{\alpha_n}| > \h \},
\end{equation*}
the process $\xi_{n\h} \equiv B(\alpha_n)$ and the filtration
$\sF_{n\h} \equiv \sG_{\alpha_n}$, so that $(\xi_t, \sF_t)_{
t \in \h\Z^+}$ is a symmetric simple random walk.  As before, 
define 
$S_t \equiv \sup_{s\leq t} \xi_s$,  $I_t \equiv \inf_{s\leq t} \xi_s$
for $t \in \h\Z^+$.

The construction  borrows the technique of \cite{R1},
where we firstly modify the given law $m$ so that the 
conditional distribution of $X_T$ given $\{
S_T=b, I_T = -a, \, \sigma_T = s \}$ is a unit mass on the
expected value $m[ X_T \,|\,  S_T=b, I_T = -a, \, \sigma_T = s\, ]$.
If we can construct a martingale with this degenerate conditional
law, then we can build the required distribution of $X_T$
given $\{S_T=b, I_T = -a, \, \sigma_T = s \}$ by Skorokhod embedding
in a Brownian motion. So we may and shall suppose
 that\footnote{There is no reason why $v$ need be a multiple
 of $h$, but this does not matter; if $s = +$, say, we shall use
 the Brownian motion living in the original probability space,
 starting at $b$ and run 
  until it first hits either the 
 upper barrier $b+\h$ or the lower barrier, which will be 
 {\em randomized}, taking value $v_+$ with
 suitably-chosen  probability $\theta$, otherwise taking value $-a-\h$.}
\begin{equation}
m[X_T = v \, |\,  S_T=b, I_T = -a, \, \sigma_T = s\, ] = 1,
\label{eq228}
\end{equation}
where $v = m[ X_T \,|\,  S_T=b, I_T = -a, \, \sigma_T = s\, ]$.

The construction is sequential, and the proof that it succeeds is
inductive. Let $\tau_n \equiv \inf\{ t: S_t - I_t = n\h\}$,
and set $\sigma_n = \alpha_{\tau_n}$, the corresponding stopping
time for the Brownian motion. The construction of $T$ begins
by setting $T = 0 $ if $U_0< m(S=I=0)$, otherwise $T \geq h = \tau_1$.
The sequential construction supposes\footnote{We
provide details of what happens if $S_{\tau_n} = \xi_{\tau_n}$;
the treatment of the case $I_{\tau_n} = \xi_{\tau_n}$ is
analogous.
} we have found that 
$T \geq \tau_n$, and $S_{\tau_n} = \xi_{\tau_n} = b$, 
$I_{\tau_n}=-a$. Then we place
a lower barrier $\ell \in [-a-h,b+h]$ by the recipe
\begin{eqnarray*}
\ell &=& v_+ 	\qquad \textrm{if } U_n < \theta
\\
&=& -a-h \qquad \textrm{else} 
\end{eqnarray*}
where $v_+$ is defined in terms of $m$ by \eqref{vdef}, and 
$\theta$ is defined by
\begin{equation}
\tilde p_{+0} \equiv
\frac{m(S=b, I = -a, \sigma = +1)}{\psi_+(-a,b)}
= \frac{m(S=b, I = -a, \sigma = +1)}{\mb(H_b<\infty,
I(S_b) = -a) }
= \theta \; \frac{\h}{b+\h-v_+}
\label{thetadef}
\end{equation}
with the notation of Proposition \ref{prop2}; in view of the
fact that we have assumed the necessary conditions \eqref{nc5}
and \eqref{nc6}, {\em we can assert\footnote{
We shall establish in the inductive proof that $\psi_\pm $
are non-negative.
} that $\theta$ so defined {\em is}
a probability}: $0 \leq \theta \leq 1$.  We now run the Brownian 
motion  $B$ forward from time $\sigma_n$ until it first hits
$\ell$ or $b+h$. If $\ell = v_+$ and $B$ hits $\ell$ before
$b+h$, then we will stop everything at that time, and 
declare that $X_T = v_+$; otherwise, we will reach either
$-a-h$ or $b+h$ and declare that $T \geq \tau_{n+1}$. If we 
determine that $T \geq \tau_{n+1}$, we take a further step of the
construction.  

For each $n \geq 1$, let $Q_n$ be the combined statement\footnote{
The functions $\psi_\pm$ are defined in terms of $m$ by
\eqref{m+}, \eqref{m-}, \eqref{psip}, \eqref{psim}.
}
\begin{center}
\begin{itemize}
\item[(i)] for all $a, b \in \h\Z^+$, $0 < a+b \leq n\h$ 
\begin{eqnarray}
P(H_b \leq T,\; I(H_b) = -a) &=& \psi_+(-a,b)
\label{Qnp}
\\
P(H_{-a} \leq T,\; S(H_{-a}) = b)
&=& \psi_-(-a,b)
\label{Qnm}
\end{eqnarray}
\item[(ii)] 
\begin{equation}
P(S=x, I = -y, X = z, \sigma = s)
= m(S=x, I = -y, X = z, \sigma = s)
\label{Qnii}
\end{equation}
 for all
$s \in\{-1,1\}$,  $x, y, z, \in \h\Z$, 
$x, y \geq 0$, $x+y < n\h$.
\end{itemize}
\end{center}
We shall prove by induction that $Q_n$ is true for all $n>0$,
establishing the statement first for $n=1$. We prove \eqref{Qnp}, 
leaving the analogous proof of \eqref{Qnm} to the diligent reader.
Taking $b = 0, \; a = \h$, 
\eqref{Qnp} says that 
\[
P(H_0 \leq T, \; I(H_0) = -\h) = \psi_+(-\h,0),
\]
and both sides are readily seen to be equal to zero; taking $b=h,\;
a = 0$, 
\eqref{Qnp} says that 
\begin{eqnarray*}
P(H_\h \leq T, \; I(H_\h) = 0) &=& \psi_+(0,\h)
\\
&=&  \varphi(\h,-\h) - \varphi(\h,0)
\\
&=&
\frac{ \h - m(\h+X; S<\h, I>-\h) }{2\h}\;  - 0 
\\
&=& \half \bigl[ \, 1 - m(S=X=I=0) \, ]
\end{eqnarray*}
which is clearly true, because if the construction does not
stop immediately at time 0 (an event of probability $m(I=X=S=0)$)
then with equal probability the process steps at time 1 to 
$\pm\h$.  The second statement \eqref{Qnii} holds because
we have constructed the probability of $I=X=S=0$ correctly.

\vskip 0.2 in
Now suppose that $Q_k$ has been proved to hold for $k \leq n$; we have
to prove \eqref{Qnp}, \eqref{Qnm} and \eqref{Qnii} for $n+1$.
To prove \eqref{Qnii}, suppose that $x, \; y \in \h\Z^+$ and 
$x+y = n\h$. By construction, the random walk will be stopped before
the range $S-I$ increases to $(n+1)\h$ if and only if the barrier
$\ell$ happens to be positioned at $v_+$ {\em and} that barrier
is hit before the Brownian motion rises to $b+h$. 
Conditional on the event $B_+ =
\{ T \geq \tau_n, \; S_{\tau_n} = \xi_{\tau_n} = b,
\; I_{\tau_n} = -a \}$, 
the probability 
of that joint event is 
\begin{equation}
    \theta \times \frac{h}{b+h-v_+}.
\end{equation}
By the inductive hypothesis \eqref{Qnp} we have that the probability 
of the conditioning event $B_+$ is $\psi_+(-a,b)$; so from the 
definition \eqref{thetadef} of $\theta$ we learn that 
\[
P(S_T = b, \; I_T = -a, \;  \sigma = +1 )
=m(S=b, I=-a, \; \sigma = +1).
\]
Given that this event happens, the conditional distribution of 
$X_T$ is correct, by the Skorohod embedding construction of $X_T$
with mean $v_+$.  Therefore \eqref{Qnii} has been proven for any
$x, \; y \in \h\Z$ with $x+y = n\h$, and for any $z \in \h\Z$, 
$s \in \{ -1, 1\}$.

\vskip 0.15 in
It remains to prove assertion (i) of $Q_{n+1}$, and for this we recall 
some of the notation of the proof of Proposition \ref{prop2}. For
$a, \; b \in \h\Z^+$, $a+b = n\h$, we  
write
\begin{eqnarray*}
p_+ &=& P(B_+ ) \equiv P(H_b \leq T, \; I(H_b) = -a),
\\
p_- &=& P(B_-) \equiv P(H_{-a} \leq T, \; S(H_{-a}) = b)
\end{eqnarray*}
which in view of the truth of $Q_n$ we know are equal to 
$\psi_+(-a,b)$ and $\psi_-(-a,b)$ respectively.
If we now define
\begin{eqnarray*}
p_{++} &=& P(B_+, H_{b+\h}  \leq T \wedge H_{-a-\h} )
\\
p_{+-} &=& P(B_+, H_{-a-\h}  \leq T \wedge H_{b+\h} )
\\
p_{+0} &=& P(B_+, T <\tau_{n+1})
\\
p_{-+} &=& P(B_-, H_{b+\h}  \leq T \wedge H_{-a-\h} )
\\
p_{--} &=& P(B_-, H_{-a-\h}  \leq T \wedge H_{b+\h} )
\\
p_{-0} &=& P(B_-, T < \tau_{n+1})
\end{eqnarray*}
then by exactly the same Optional Sampling argument which led to
\eqref{p++}, \eqref{p+-}, we conclude that
\begin{eqnarray}
p_{++} &=& \frac{(b+a+\h)p_+ - (a+\h+v_+)\, p_{+0}}{b+a+2\h}
\label{ppp}
\\
p_{+-} &=& \frac{\h p_+ - (b+\h-v_+)\,  p_{+0}}{b+a+2\h}
\label{ppm}
\\
p_{-+} &=& \frac{\h p_- - (a+\h+v_-) p_{-0}}{a+b+2\h}
\label{pmp}
\\
p_{--} &=& \frac{(a+b+\h)p_- - (b+\h-v_-)p_{-0}}{a+b+2\h}
\label{pmm}
\end{eqnarray}
and now the task is to prove (after cross-multiplying
by $a+b+2\h$) that 
\begin{equation}
(a+b+2\h)\{\,p_{++} + p_{-+}\, \}
= (a+b+2\h)\psi_+(-a, b+\h),
\label{toprove}
\end{equation}
and the minus analogue, which is just the same argument
{\it mutatis mutandis.}
Firstly we develop the left-hand side using \eqref{ppp}, 
\eqref{ppm} and their analogues for $B_-$ to obtain
\begin{eqnarray*}
LHS &=& (a+b+\h)\psi_+(-a,b) - (a+\h+v_+)p_{0+}
+\h\psi_-(-a,b) - (a+\h+v_-) p_{-0}
\\
&=& (a+b+\h)\{ \, \vp(b,-a-\h) - \vp(b,-a)\, \}
+ \h\{ \, \vp(-a,b+\h) - \vp(-a,b)\, \}
\\
&&\qquad\qquad\qquad
-(a+\h) m(S=b,I=-a) - m(X; S = b, I= -a)
\\
&=& a+\h-m(a+\h+X; S<b, I>-a-\h) - \{\, a - m(a+X; S<b, I>-a)\,\}
\\
&&\qquad\qquad
- \h(\vp(b-a) + \vp(-a,b)) + \h\vp(-a,b+\h) 
-m(a+\h+X; S=b, I = -a)
\\
&=& \h-m(a+\h+X; S<b, I>-a-\h) +m(a+X; S<b, I>-a)
\\
&&\qquad\qquad
-\h\{1 - m(S<b,I>-a)\}+ \h\vp(-a,b+\h)-m(a+\h+X; S=b, I = -a)
\\
&=& -m(a+\h+X; S<b, I>-a-\h) +m(a+\h+X; S<b, I>-a)
\\
&&\qquad\qquad
-m(a+\h+X; S=b, I = -a) + \h\vp(-a,b+\h)
\\
&=& -m(a+\h+X: (A_2\cup A_3) \backslash A_1) + \h\vp(-a,b+\h)
\end{eqnarray*}
where $A_1 = \{ S<b, I > -a\}$, $A_2 = \{
S<b, I>-a-\h\}$ and $A_3 = \{ S=b, I = -a\}$. Noticing that  
$A_1 \subseteq A_2$ and $A_3$ is disjoint from $A_1$,
the region of integration is
\[
(A_2\cup A_3) \backslash A_1= \{ S<b, I = -a\} \cup A_3
= \{ S \leq b, I = -a\} 
= \{ S < b+\h, I = -a\}.
\]  
Hence the left-hand side is equal to 
\begin{equation}
LHS = -m(a+\h+X; S<b+\h, I = -a) + \h\vp(-a,b+\h).
\label{LHS}
\end{equation}

\medskip\noindent
Turning now to the right-hand side of \eqref{toprove}, we have
\begin{eqnarray}
RHS &=& (a+b+2\h)\{ \, \vp(b+\h,-a-\h) - \vp(b+\h,-a)\,\}
\nonumber
\\
&=& a+\h-m(a+\h+X:S<b+\h,I>-a-\h) -\h\vp(b+\h,-a)
\nonumber
\\
&&\qquad\qquad\qquad
- \{ \, a - m(a+X:S<b+\h,I>-a)\,\} 
\nonumber
\\
&=& \h - m(a+\h+X:S<b+\h,I>-a-\h)+m(a+\h+X;S<b+\h,I>-a)
\nonumber
\\
&&\qquad\qquad
-\h m(S<b+\h,I>-a) - \h\vp(b+\h,-a)
\nonumber
\\
&=&\h\{\, 1-m(S<b+\h,I>-a) - \vp(b+\h,-a)\, \}
\nonumber
\\
&&\qquad\qquad\qquad
- m(a+\h+X;S<b+\h, I = -a).
\label{RHS}
\end{eqnarray}
Comparing \eqref{LHS} and \eqref{RHS}, we see that we have to prove
\begin{equation}
\vp(b+\h,-a) + \vp(-a,b+\h) = 1- m(S<b+\h,I>-a),
\end{equation}
which is evidently true from the definition \eqref{m+}, 
\eqref{m-} of $\vp$.

\hfill $\square$

    \section{Hedging.}\label{hedge}
Theorem \ref{thm1} provides us with necessary and sufficient
conditions for a measure $m$ on $\sX$ to be consistent. In
principle, this allows us to construct extremal martingales, and 
robust hedges for derivatives.  

Let us firstly see how this works in the context of the joint law
of $(S,X)$ studied in \cite{R1}. We begin by recalling some of 
the results of that paper.
We let $X_t = B_{t \wedge T}$ be
a Brownian motion stopped as an almost-surely finite stopping
time $T$, with $S_t = \sup_{u \leq t} X_u$, and with
$S \equiv S_\infty$, $X \equiv X_\infty$.  With this terminology,
Theoren 3.1 of \cite{R1} says the following.

\begin{theorem}
The probability measure $\mu$ on $\R^+ \times \R^+$ is the 
joint law of $(S, S-X)$ for some almost-surely finite
stopping time $T$ if and only if
\begin{equation}
\biggl(  \iint_{ (t,\infty) \times \R^+} \mu(ds,dy)
\biggr) dt  \geq \int_{(0,\infty)} y \; \mu(dt, dy).
\label{R1_3.1}
\end{equation}
If $(X_t)_{t\geq 0}$ is also uniformly integrable, then
inequality \eqref{R1_3.1} holds with equality:
\begin{equation}
\biggl(  \iint_{ (t,\infty) \times \R^+} \mu(ds,dy)
\biggr) dt  = \int_{(0,\infty)} y \; \mu(dt, dy).
\label{R1_3.2}
\end{equation}
Finally, if \eqref{R1_3.2} holds, and if $X \in L^1$,
\begin{equation}
\iint |t-y| \; \mu(dt,dy) < \infty,
\label{X_in_L1}
\end{equation}
 then $\mu$ is the joint law of $(S, S-X)$ for a
uniformly integrable martingale $(X_t)_{t\geq 0}$.
\end{theorem}

\medskip
\noindent
{\sc Proof.} See \cite{R1}. The final assertion is not 
in \cite{R1}, but can easily be deduced. In view of the first
assertion, there is some stopping time $T<\infty$ such that 
$\mu$ is the joint law of $(S,S-X)$. By multiplying \eqref{R1_3.2}
by some non-negative test function $\varphi$ and integrating 
with respect to $t$ we discover that
\begin{equation}
\mu(\Phi) = \mu( \, (S-X)\varphi(S) \,)
\label{eq34}
\end{equation}
where $\Phi(t) = \int_0^t \, \varphi(y) \; dy$. Taking
$\varphi(x) = I_{\{ x > b \}}$ for some $b\geq 0$ we find that
\begin{equation}
b \mu(S>b)  = \mu( X: S>b).
\label{eq35}
\end{equation}
Using the fact that $X \in L^1$, we can let $b \uparrow \infty$ 
in \eqref{eq35} to prove that $\lim_{b \uparrow\infty}
b \mu(S>b) =0$. Lemma 2.3 of \cite{R1} gives the result.

\hfill$\square$

\medskip\noindent
{\sc Remark.} Standard monotone class arguments show that 
\eqref{R1_3.1} is equivalent to the statement that 
\begin{equation}
 \mu(\Phi) \geq \mu( \, (S-X)\varphi(S) \,)
 \label{suff1}
\end{equation}
for all non-negative test functions, which again is equivalent
to the statement that 
\begin{equation}
 b \mu(S>b)  \geq \mu( X: S>b)
 \label{suff2}
\end{equation}
for all $b \geq 0$.   Likewise, \eqref{R1_3.2} is equivalent 
to \eqref{eq34} for all non-negative test functions $\varphi$,
which again is equivalent to the statement \eqref{eq35}:
\begin{equation}
\mu( \, X-b : S>b ) = 0\qquad\forall b \geq  0.
\label{suff3}
\end{equation}

\vskip 0.1 in
An important and typical\footnote{
The papers  Hobson\cite{Hobson1}, ... give examples of this
kind.
} use of this would be to try to find
an {\em extremal} martingale, which would in turn lead to a
maximum possible
derivative  price and a robust hedging strategy. So, for
example, suppose that we observe call  option prices $C(K)$ for every
strike $K$ at a common fixed expiry time\footnote{
Let us suppose that the expiry is 1.} for some (discounted)
asset, and suppose that the asset has continuous paths
 $(X_t)_{0 \leq t\leq 1}$, and 
 is a uniformly-integrable martingale in the 
pricing measure. 

 Suppose now that we are given some
 derivative whose payoff at time 1 is $G(S_1, X_1)$, where
 $S_1 = \sup_{0 \leq t \leq 1} X_t$; {\em what is the most
 expensive the time-0 price of this derivative can be?}
 
 The time-0 price of the derivative is given by
 \begin{equation}
\iint G(s, x)  \; q(ds,dx)
\label{obj1}
\end{equation}
where $q$ is the joint law\footnote{
As before, when the time subscript of a process is omitted, we
understand it to be 1.
} of $(S, X)$. Now provided the law $q$ satisfies the conditions
\begin{equation}
\iint (x-K)^+ \; q(ds,dx) = C(K) \qquad\forall K
\label{cons1}
\end{equation}
and  (see \eqref{suff3})
\begin{equation}
\iint_{s>b} (x-b) \; q(ds,dx) = 0 \qquad\forall b>0
\label{cons2}
\end{equation}
then $q$ is the  joint distribution of $(S, X)$ for {\em some }
continuous martingale whose law at time 1 agrees with the data
contained in the call prices.  The problem of finding the most
expensive time-0 price is therefore the problem of maximizing the
{\em linear} objective \eqref{obj1} over non-negative probability
measures $q$ subject to the {\em linear}
 constraints \eqref{cons1} and 
\eqref{cons2}.  
Writing the problem in Lagrangian form\footnote{
This linear programming approach to the problem is also used
in \cite{DOR}.
}, we seek
\begin{eqnarray}
L(\alpha,\eta,\lambda) &=& \sup_{ q \geq 0} \biggl[\;
\iint \bigl\lbrace\, G(s,x) - \alpha - \int (x-K)^+ \;
\eta(dK) + \int_0^\infty (x-b)I_{\{ s > b\}} \; \lambda(db)
\bigr\rbrace \; q(ds,dx)
\nonumber
\\
&&\qquad\qquad\qquad + \alpha + \int C(K) \; \eta(dK) \; \biggr].
\label{lagr1}
\end{eqnarray}
From standard linear programming results, we would expect that
for dual feasibility we must have 
\begin{equation}
G(s,x) \leq \alpha + \int (x-K)^+ \;
\eta(dK) - \int_0^\infty (x-b)I_{\{ s > b\}} \; \lambda(db)
\label{robust_hedge}
\end{equation}
everywhere, with equality everywhere that the optimal $q$ places
mass; and that the dual problem will be
\begin{equation}
\inf \biggl[ \; 
\alpha + \int C(K) \; \eta(dK)  \; \biggr]
\label{dualLP}
\end{equation}
over $(\alpha,\eta,\lambda) $ satisfying \eqref{robust_hedge}.
These equations have a simple and beautiful interpretation. 
The dual-feasibility relation \eqref{robust_hedge} expresses
a {\em robust hedge}; if we hold $\alpha$ in cash, $\eta(dK)$
calls of strike $K$, and {\em sell forward $\lambda(db)$ units of
the underlying when $S$ reaches the level $b$}, then we generate
a contingent claim at the terminal time which will always 
dominate the claim $G$ which we have to pay out.   The dual
form of the linear program \eqref{dualLP} says that the cost of
constructing such a hedge, which is of course 
$\alpha + \int C(K) \; \eta(dK)$, must be minimized.

The  primal problem seeks to find the most expensive that the
derivative $G(S,X)$ can be, given the market prices $C(K)$; and the
dual problem seeks the cheapest super-replicating hedge.  The
characterization \eqref{suff3} of the possible joint laws of $(S,X)$
{\em tells us what the form of the hedge \eqref{robust_hedge} must be.}

\vskip 0.2 in

Our goal now is to try to use  Theorem \ref{thm1} to similarly
bound the price of, and to super-replicate, contingent claims
which depend on the maximum, terminal value, {\em minimum,
and direction of the final excursion} for a stopped symmetric
simple random walk. To understand how this is to be done, we 
focus on the `plus' versions of the necessary and sufficient
conditions \eqref{nc5}. We shall also suppose that the 
martingale $X$ is {\em uniformly integrable}, to avoid 
having to bother about side issues.

The condition \eqref{nc5} can be restated in terms of the measure
$m$ as 
\begin{eqnarray}
m(b+h-X: S=b, I=-a, \sigma=+1) &\leq & h \psi_+(-a,b)
\label{HE1}
\\
 &=& h \{ \,\varphi(b,-a-h) - \varphi(b,-a) \, \}
 \nonumber
\end{eqnarray}
in the notation of Section \ref{S1}.  From the definition
\eqref{m+} of $\varphi(b,-a)$, from the fact that $m(X)=0$, and
the Optional Sampling Theorem result that $m(a+X:I \leq -a)=0$,
we have 
\begin{eqnarray*}
(a+b)\varphi(b,-a) &=& a - m(a+X: S <b, I > -a)
\\
&=& m(a+X: \textrm{ $S\geq b$ or $I \leq -a$} )
\\
&=& m(a+X: S \geq b, I > -a)
\\
&=& (a+b) m(S \geq b, I > -a) - m(b-X: S \geq b, I > -a).
\end{eqnarray*}
Thus the inequality \eqref{HE1} may be re-expressed  after
some simple rearrangement as
\begin{eqnarray}
0 &\leq & h m(S\geq b, I = -a) - \frac{h}{a+b+h} \; m(b-X:
S\geq b, I > -a-h) +
\nonumber
\\
&&\qquad + \frac{h}{a+b} \; m(b-X: S \geq b, I > -a)
- m(b+h-X: S=b, I=-a, \sigma=+1).
\nonumber
\end{eqnarray}
This inequality for all $a,\, b \in h\Z^+$ not both zero, together
with the `minus' analogues, is necessary and sufficient for a 
probability measure $m$ to be the joint law of $(I,X,S,\sigma)$.
Just as we did  at \eqref{lagr1}
for derivatives depending only on $(X,S)$,
we can construct the Lagrangian for this problem, which would 
give us terms of the form
\begin{eqnarray}
\lambda^+_{ab} \;(Z-w) &\equiv &
\lambda^+_{ab}\biggl[ \;
h I_{\{S \geq b, I = -a\}} - \frac{h}{a+b+h} \;(b-X)
I_{\{ S \geq b, I > -a-h\}} +
\nonumber
\\
&&  
+ \frac{h}{a+b} \; (b-X) I_{\{ S \geq b, I > -a \}}
 - (b+h-X)I_{\{ S=b, I=-a, \sigma=+1 \}}  -w
 \;\biggr],
 \label{Zdef}
\end{eqnarray}
where $w\geq 0$ is a non-negative slack variable to handle the
inequality constraint.  Dual feasibility will therefore require
that $\lambda^+_{ab}\geq 0$, and at optimality we will have 
the complementary slackness condition $\lambda^+_{ab}\, w = 0$.

In the situation of derivatives depending only on $(X,S)$, we
had  terms of the form $\lambda_a (X-a) I_{\{S>a\}}$, 
which were interpreted as forward purchase of the underlying
asset when the supremum process reaches a new level. This
forward purchase interpretation determines a hedging strategy
which {\em can be implemented in an adapted fashion}. However, 
it is very far from clear that the random variable $Z$ defined
at \eqref{Zdef} can be realized by some adapted 
trading strategy. For example, the term involving
$(b-X) I_{\{ S \geq b, I > -a \}}$ could be interpreted as a forward
sale of the underlying when the price first gets to $b$; but
this trade should only be put on if $I>-a$, and it is not known
at time $H_b$ whether or not the ultimate infimum $I$ will be
greater than $-a$ or not.

Nevertheless, we can specify an adapted trading strategy which
will subreplicate the random variable $Z$, as follows.  We construct
a random variable $Y$ which is the final value of the adapted
hedging strategy made up of three component positions:
\begin{enumerate}
\item At $H_b$, buy forward $h/(a+b+h)$ units of the underlying
if $I(H_b)>-a-h$, and come out of the position at time $H_{-a-h}$;
\item At $H_b$, buy forward $-h/(a+b)$ units of the underlying
if $I(H_b)>-a$, and come out of the position at time $H_{-a}$;
\item At $H_b$, buy forward 1 unit of the underlying
if $I(H_b)=-a$, and come out of the position at time $H_{b+h}
\wedge H_{-a-h}$.
\end{enumerate}
Now clearly the random variable
\begin{eqnarray}
Z &\equiv & h I_{\{S \geq b, I = -a\}} - \frac{h}{a+b+h} \;(b-X)
I_{\{ S \geq b, I > -a-h\}} +
\nonumber
\\
&&  
+ \frac{h}{a+b} \; (b-X) I_{\{ S \geq b, I > -a \}}
 - (b+h-X)I_{\{ S=b, I=-a, \sigma=+1 \}}
 \label{Zdef2}
\end{eqnarray}
will be zero if $S < b$ or if $I\leq -a-h$, so to understand $Z$
we may suppose that $H_b<\infty = H_{-a-h}$. 

But before we narrow 
our attention down to the event $\{ H_b<\infty = H_{-a-h}\}$, we
should consider what happens off that event to $Y$. If $H_b = 
\infty$, then none of the component positions of $Y$ is ever
entered, so $Y=0$ in that case.  If $H_b<\infty$ and $H_{-a-h}
< \infty$, then we have three cases to consider:
\begin{itemize}
\item[(i)] When $I(H_b)>-a$, the strategy enters positions 1 and 2
at time $H_b$, and closes out both when the infimum falls to $-a$ and
then to $-a-h$;  position 1 loses $h$, position 2 gains $h$, so altogether
$Y=0$;
\item[(ii)] When $I(H_b)=-a$, the strategy enters positions 1 and 3.
If $H_{b+h}< H_{-a-h}$, then position 3 makes a gain of $h$ when it is 
closed out, but position 1 makes a loss of $h$ when it is closed out, 
so overall zero gain. On the other hand, if $H_{-a-h}<H_{b+h}$, then 
position 1 makes a loss of $h$ when it is closed out, and 
position 3 makes a loss of $(a+b+h)$ when it is closed out, so overall
$Y = -(a+b+h)-h<0$, and as we shall
subsequently see, {\em  this is the only situation in which $Y$ is 
strictly less than $Z$};
\item[(iii)] When $I(H_b) \leq -a-h$, none of the positions is entered, 
and $Y=0$.
\end{itemize}

We now have to compare the values of $Z$ and $Y$ on the event 
$\{ H_b<\infty = H_{-a-h}\}$, breaking the comparison down into seven
cases as presented in the following table. In the first two rows, we 
see what happens if $I>-a$, and in the remaining rows, we are
considering situations where $I=-a$. The reader is invited to check
through each of the entries of the table, and confirm the findings
reported there. The only entry that requires comment is the penultimate
row, in the column for $Z$. In this row, we are in the situation
where $S=b$ and $I=-a$, so we get a contribution to $Z$ from the
first term in \eqref{Zdef2}, and from the second term, none from 
the third term, and {\em none from the fourth term}, because if 
$H_b<H_{-a}<H_{b+h} = \infty$ it must be that {\em the signature
$\sigma$ is $-1$\, !}  What we see from the table is that in every
case the value of $Z$ is equal to the value of $Y$.

\bigbreak
\begin{center}
\begin{tabular}{|c| c |c|}
\hline 
$H_{-a-h}=\infty$ & $Z$ & $Y$ \\ 
\hline 
\hline
$H_b<H_{b+h}<\infty=H_{-a}$ & $\frac{h(b-X)}{a+b}-\frac{h(b-X)}{a+b+h}$ & 
$\frac{h(X-b)}{a+b+h}-\frac{h(X-b)}{a+b} $\\ 
\hline 
$H_b<H_{b+h}=\infty = H_{-a}$ & $\frac{h(b-X)}{a+b}-\frac{h(b-X)}{a+b+h}$ 
& $\frac{h(X-b)}{a+b+h}-\frac{h(X-b)}{a+b} $ \\ 
\hline 
$H_{-a}<H_b<H_{b+h}<\infty$ & $h-\frac{h(b-X)}{a+b+h}$ & 
$ \frac{h(X-b)}{a+b+h} +h $ \\ 
\hline 
$H_{-a}<H_b<H_{b+h}=\infty $ & $h - \frac{h(b-X)}{a+b+h}+(X-b-h)$ & 
$\frac{h(X-b)}{a+b+h} +X-b $ \\ 
\hline 
$H_b<H_{-a}<H_{b+h}<\infty$ & $h - \frac{h(b-X)}{a+b+h}$ 
& $\frac{h(X-b)}{a+b+h} +h  $ \\ 
\hline 
$H_b<H_{-a}<H_{b+h}=\infty$ & $h-\frac{h(b-X)}{a+b+h}$ 
& $\frac{h(X-b)}{a+b+h} +h  $ \\ 
\hline 
$H_b<H_{b+h}<H_{-a}<\infty$ & $h - \frac{h(b-X)}{a+b+h}$ 
&  $\frac{h(X-b)}{a+b+h} +h  $ \\ 
\hline 
\end{tabular} 
\end{center}

Thus we may conclude that $Y \leq Z$ in all instances, and the only
situation in which the inequality is strict is when $H_{-a} < H_b
< H_{-a-h} < H_{b+h}$.

Now we explain how these observations lead to a super-replicating
hedging strategy. For this, let us denote by $Z^+_{ab}$ then 
random variable we have been calling $Z$ up til now; this is because
in the Lagrangian we have to consider such random variables (and their
`minus' analogues) for all $a, \, b \in h\Z^+$ not both zero.  Suppose
that we have some derivative $G(I,X,S,\sigma)$ whose price we wish to 
maximize subject to the distribution of $X$ matching call price
data, just as we did for derivatives depending only on $(X,S)$ in the
first part of our discussion in this Section.  We would find ourselves
with a Lagrangian form similar to \eqref{lagr1}:
\begin{eqnarray}
L(\alpha,\lambda,\eta) &=& 
\sup_{m\geq 0}\biggl[\;
\int \bigl\lbrace\;
G(I,X,S,\sigma) - \alpha - \int (X-K)^+ \eta(dK)
+ 
\nonumber\\
&& \qquad  
+ \sum_{a,b,\pm} \lambda^\pm_{ab}(Z^{\pm}_{ab} - w^\pm_{ab})
\;\bigr\rbrace \; dm(I,X,S,\sigma)
+ \alpha + \int C(K) \; \eta(dK)
\;\bigr\rbrace
\;\biggr]
\end{eqnarray}
with obvious notation.  Now dual feasibility imposes the condition
\begin{eqnarray}
G(I,X,S,\sigma) & \leq &
\alpha + \int (X-K)^+ \eta(dK)-\sum_{a,b,\pm} \lambda^\pm_{ab}\,Z^{\pm}_{ab}
\label{Zhedge}
\\
&\leq &\alpha + \int (X-K)^+ \eta(dK)-\sum_{a,b,\pm} \lambda^\pm_{ab}\,Y^{\pm}_{ab}
\label{Yhedge}
\end{eqnarray}
in another obvious notation.  The interpretation of \eqref{Yhedge} is 
that {\em the derivative $G$ is super-replicated by the adaptively-realizable
hedge given by a position in calls and a position in the $Y$-hedges.}

At optimality, complementary slackness tells us that if
$\lambda^+_{ab}>0$ then $w^+_{ab}=0$, and therefore the
inequality \eqref{HE1} must hold with equality. Tracing this back to 
the condition \eqref{nc5}, and its derivation from \eqref{p+-}, we 
find that equality in \eqref{HE1} is equivalent to the statement 
that $\tilde p_{+-}=0$.  What this means is that on the event
$\{ H_{-a}<H_b <H_{-a-h}\}$ we {\em cannot have $H_{-a-h}<H_{b+h}$},
and as we saw, this was {\em the only situation where $Y<Z$}. We 
may therefore conclude that for the optimal $m^*$, not only will
\eqref{Zhedge} hold with equality $m^*$-a.e., but also \eqref{Yhedge}
will hold with equality $m^*$-a.e..  In other words, if the joint
law $m$ is the optimal joint law, the hedging strategy expressed
by \eqref{Yhedge} is a perfect replication of the contingent claim - there
is no slack.


\pagebreak
\bibliography{IXSbib}

\begin{thebibliography}{1}

\bibitem{BHR1}
H.~Brown, D.~Hobson, and L.~C.~G. Rogers.
\newblock The maximum maximum of a martingale constrained by an intermediate
  law.
\newblock {\em Probability {T}heory and {R}elated {F}ields}, 119:558--578,
  2001.

\bibitem{BHR2}
H.~Brown, D.~Hobson, and L.~C.~G. Rogers.
\newblock Robust hedging of barrier options.
\newblock {\em Mathematical Finance}, 11:285--314, 2001.

\bibitem{CO}
A.~M.~G. Cox and J.~Obloj.
\newblock Robust pricing and hedging of double no-touch options.
\newblock {\em Finance and Stochastics}, 15:573--605, 2011.

\bibitem{DOR}
M.~H.~A. Davis, J.~Obloj, and V.~Raval.
\newblock Arbitrage bounds for prices of weighted variance swaps.
\newblock {\em Mathematical Finance}, ??:???--???, 20??

\bibitem{Feller}
W.~Feller.
\newblock {\em An {I}ntroduction to {P}robability {T}heory and its
  {A}pplications}, volume~1.
\newblock Wiley, New York, third edition, 1968.

\bibitem{Hobson1}
D.~G. Hobson.
\newblock The maximum maximum of a martingale.
\newblock In J.~Az{\'e}ma and M.~Yor, editors, {\em S{\'e}minaire de
  {P}robabilit{\'e}s}, volume XXXII, pages 25--263. Springer, 1998.

\bibitem{R1}
L.~C.~G. Rogers.
\newblock The joint law of the maximum and the terminal value of a martingale.
\newblock {\em Probability Theory and Related Fields}, 95:451--466, 1993.

\end{thebibliography}
\bibliographystyle{plain}

\end{document}